\documentclass{article}

\newtheorem{df}{ \sc Definition}[section]

\newtheorem{cl}[df]{ \sc Claim}
\newtheorem{thrm}[df]{ \sc Theorem}

\newtheorem{re}[df]{ \it Remark}

\newtheorem{cj}[df]{ \sc Conjecture}
\newtheorem{constr}[df]{ \sc Construction}

\def\dim{{\rm dim}}

\def\det{{\rm det}}

\def\Pic{\rm Pic}
\def\epi#1{\;\smash{\mathop{\hbox to 20pt{\rightarrowfill}\hskip
-15pt\rightarrow}\limits^{#1\,}}\;}
\def\mono{\lhook\joinrel\relbar\joinrel\rightarrow}
\def\exseq#1#2#3{0\rightarrow #1 \rightarrow #2 \rightarrow #3 \rightarrow 0}
\def\O{{\mathcal O}}
\def\Ac{{\mathcal A}}
\def\Bc{{\mathcal B}}
\def\Mc{{\mathcal M}}
\def\H{{\mathcal H}}
\def\C{{\mathbb C}}
\def\P{{\mathbb P}}
\def\Z{{\mathbb Z}}
\def\A{{\mathbb A}}

\def\Q{{\mathbb Q}}

\def\Proof{\noindent\hskip4pt{\it Proof}.\ }
\def\qed{\hfill$\Box$\vskip10pt}                                                           

\input{amssymb.sty}
\usepackage{amscd}  
\usepackage{amsmath}                                                  
\begin{document}
\hsize=5.3in
\author{Nicolae Manolache
} 
\title{Cohen - Macaulay Nilpotent Schemes}
\date{}
\maketitle

\begin{abstract} 
We present here a short and partial survey about the construction and classification
of the Cohen-Macaulay scheme structures on a smooth variety as support, or
on a union of smooth varieties. We present, in the chronological order, the
results of {\bf Fossum} (cf. \cite{Fo}), {\bf Ferrand} (cf. \cite{Fe}), {\bf B\u anic\u a-Forster}
(cf. \cite{BF1}, \cite{BF2}) and those of \cite{M2}, \cite{M5}. We note that in the last ten 
years a new terminology emerged,
as for instance {\bf ribbon} for a double structure on a scheme (usually of pure dimension 1)
considered as an abstract scheme, i.e. not embedded, (there is a notion of {\em ribbon}
also in the algebraic topology) and {\bf rope} which is
a multiple structure $Y$ on a scheme $X$ whose associated graded algebra has only 
one non-trivial component (i.e. any of the canonical filtrations is simply $Y\supset X$).
The introducing of the new names created islands of research which unfortunately do not
communicate among them.

\end{abstract}

\section{Introduction}
I was asked by the organizers of this Seminar to present a survey, for nonspecialists 
in algebraic geometry, of a subject which is not intuitive in the classical sense.
I try here to realize this job, without giving too many formal definitions, and without
 complete proofs. The interested reader is invited to follow the references to fill
the gaps. The examples of how to apply the general theory given in the last section are 
more technical, although "elementary".
The multiple structures considered there are of the simplest possible kind,
being {\em primitive} (see the definition below). 
In order to see more applications and examples we invite the reader to go to the papers in the
{\em References} and to the references to those papers

In many geometric problems, degenerated cases appear naturally, even in
elementary circumstances, as, for instance,
double (plane) lines as degenerated conics, sets of points among which
some are "thicker" points (e.g. intersection of curves having a high contact
in some points), etc. 
These circumstances made one of the difficulty of the classical algebraic 
geometry.

I will remind some questions of interest in Algebraic Geometry, where nilpotent 
structures arise naturally:

(i) The concept of (algebraic) schemes introduced by A. Grothendieck
in the years 1950's as the next step to Serre's famous paper \cite{Se1} 
made possible to attach to an algebraic (projective, affine,
etc. ) set different structures, so to say "according to the equations defining it".
Standard nilpotent structures are the {\em infinitesimal neighbourhoods}.
If $X$ is an affine (or projective) scheme, defined in
the affine (projective) space by an ideal $I$, then the scheme $X^{(1)}$ defined by $I^2$
is calleded {\em the first infinitesimal neighbourhood}. Also the scheme $X^{(2)}$
defined by $I^3$ is called the second infinitesimal neighbourhood, etc. 
These are given by very many equations and the multiplicity increases dramatically. 
We don't give here the definition of the multiplicity, but in the case of smooth $X$ in a 
projective space the multiplicity of a mild (i.e. Cohen-Macaulay, see below) structure on $X$  
is the factor by which the degree of $X$ is multiplied. This factor is a positive integer.
For instance the question posed by Kronecker and Severi about the minimal number
of equations defining a smooth affine space curve meant in this new frame
to find a {\bf new structure, a fortiori nilpotent} on the curve, having as few 
equations as possible. 

(ii) The frame of schemes made possible the definition of new objects as, for instance
{\bf the Hilbert scheme}. If we consider a projective algebraic set $X \subset \P^n$ 
(zero set of some homogeneous polynomials), one can attach to it a graded ring taking the
quotient of the  polynomial ring by the ideal generated by all homogeneous
polynomials which are zero on $X$. The dimension of the component of degree $N$,
for $N>>0$ is polynomial in $N$ (consequence of the Hilbert syzygies theorem).
This polynomial is called {\em the Hilbert polynomial of $X$} and its properties are
related to the properties of $X$. For instance its degree is {\em the dimension}
of $X$. In the case of smooth projective curves, say over $\C$, the Hilbert polynomial 
has the form $ \chi (N)= dN+1-g$, where $d$ is the degree (it depends upon the projective
embedding) and $g$ is the {\em genus}, which is defined {\bf topologically}.
The set of {\bf all} projective subschemes of $\P ^n$ with a fixed Hilbert polynomial
is shown to be naturally an algebraic scheme, called {\em Hilbert scheme}. 
By a theorem of Hartshorne, the Hilbert schemes are connected. Intuitively, this means 
that one can deform a given projective subscheme of $\P^n$ to any other one with the 
same Hilbert polynomial. But the Hilbert scheme, even in the case of curves (Hilbert polynomial 
of degree $1$) has points which correspond to too wild objects, as for instance unions of 
curves and points. Even worst, the suplementary points can be {\em embedded}. To have an example, 
of what an embedded point means, consider the algebraic set $X$ defined by the 
ideal $I:=(x^2,xy)$ in the projective plane $\P^2$ with homogeneous coordinates $x,y,z$. 
As a set, this is the line $d$ given by $x=0$,
but in fact the primary decomposition $I=(x)\cap (x^2,y)$ shows  that "scheme-theoretically" 
$X=d\cup Q$, where $Q$ is the {\bf double point} $y=0, x^2=0$.
 
A natural question is :  is the locus of smooth curves connected ? (i.e. is the set of
points in the Hilbert scheme, corresponding to the smooth curves with a fixed Hilbert
polynomial connected?) An example of Gruson and Peskine (cf. \cite{GP}) shows that for
$d=9, g=10$ this locus has two irreducible disjoint components.

Hartshorne  refined this question to the following one:
If we fix a Hilbert polynomial, is the locus of the projective subschemes of $\P^3$ of 
(pure) dimension $1$ (i.e. all the irreducible components are of dimension $1$) and 
without {\bf embedded points} connected ?
In order to answer this question, one has to consider (mild) nilpotent curves.
The answer is known to be affirmative only for $d\le 4$, cf. \cite{N}, \cite{NS}.

(iii) A method to study vector bundles $E$ (the method is efficient for rank two) is to 
consider the zero-sets
of its sections (or zero-sets of a twist $E(\ell):=E\otimes \O(\ell)$, 
by a power $\O(\ell)=\O(1)^\ell$ of the tautological bundle $\O(1)$. If $\ell << 0$,
$\Gamma (E(\ell))=0$. If we take the smallest value of $\ell$ for which one has nonzero
sections, the associated zero-sets will have pure  codimension equal to the rank of $E$
and will be locally complete intersection. In general one gets a scheme structure with
nilpotents (cf. \cite{H1},\cite{M1},\cite{BM},\cite{HVdV},\cite{M4},\cite{OSz}).

(iv) One of the most challanging conjectures in the Projective Algebraic Geometry is the 
following one, due to Hartshorne, considered here in the special case of codimension $2$:
\begin{cj}\label{ci}(Hartshorne)
If $X\subset \P=\P^N,\ N\ge 6$ is a smooth projective variety of codimension $2$,
then $X$ is (globally) complete intersection.
\end{cj}

Equivalent to this is the following

\begin{cj}\label{vb} (Rank $2$ bundle conjecture)
If $E$ is a rank $2$ vector bundle on $\P=\P^N,\ N\ge 6$, then $E$ splits.
\end{cj}

Surprisingly enough, equivalent to these two is the following one:

\begin{cj}\label{stm}
If $Y$ is a multiplicity $3$ structure (see the definition bellow) on a linear
subspace $\P=\P^N$ of $\P^{N+\delta}$, contained in the first infinetisimal neighbourhood
of $\P$, then there exists a double structure $Z$ on $\P$ contained in $Y$.\end{cj}

For other equivalent variants see \cite{V2}.

\section{Fossum - Ferrand doubling}

\subsection{The algebraic case (Fossum)}
For us all the rings $R$ will be local commutative noetherian algebras over an algebraically 
closed field $k$. In this case a local ring $R$ is called Cohen-Macaulay iff the maximal 
ideal $m_R$ of $R$ contains $n=\hbox{dimension}(R)$ elements $x_1,\ldots ,x_n$ which 
"behave as {\em indeterminates}". Here {\em a sequence of elements $x_1,\ldots ,x_m$ 
behaves as indeterminates"} means
 
(a) {\em the homomorphism sending $X_i$ in $x_i$ is a faithfully flat injection \\
$k[X_1,\ldots ,X_m]\rightarrow R$ from the ring of polynomials in $m$ indeterminates 
into $R$}

or, equivalently:

(b) {\em $x_1$ is not a 
zero-divizor in $R$, $x_2$ is not a zero-divizor in $R/x_1R$, \ldots , $x_m$ is not a 
zero-divizor in $R/(x_1R+\ldots +x_{m-1}R)$}. 

A sequence of such elements $x_1,\ldots ,x_m \in m_R$, $m\le n$, is 
called {\em regular}.

Over a ring $R$, the functor of taking the dual with respect to $R$ (i.e. $Hom(?,R)$)
does not have the properties which we would require for a good duality. The Cohen-Macaulay
rings have a {\em dualizing module} (called also {\em canonical
module}) $E$.  A canonical module $E$ over a commutative noetherian local ring $R$ is a 
finitely generated module with the properties:

(a) the canonical morphism $R\rightarrow Hom_R(E,E)$ is bijective

(b) for all $i>0$, $Ext_R^i(E,E)=0$

(c) E has finite injective dimension.

One shows also that, conversely, if $R$ has a canonical module, then $R$ is Cohen-Macaulay.
Of special interest are the rings which are canonical modules over themselves. They are
called {\em Gorenstein rings}. Regular rings and complete intersection rings are Gorenstein.

One shows that a quotient $S$ of a Gorenstein ring $R$ is Cohen-Macaulay iff $Ext^i_R(S,R)=0$
for all $i\neq d:=\hbox{dim}R-\hbox{dim}S$. In this case the canonical module of $S$ is
$Ext^d_R(S,R)$.

Reiten proved in \cite{R} the following: 
\begin{thrm}\label{Reiten} Let $S$ be a Cohen-Macaulay ring and $M$ a canonical $S$-module.
Then the trivial extension $M\rtimes S$ of $S$ by $M$ is Gorenstein.{\em 
The trivial extension } means the ring structure on $M\times S$ given by $(s_1,m_1)(s_2,m_2)=
(s_1m_2+s_2m_1, s_1s_2)$.
\end{thrm}
\qed

More generally, if $S$ is a ring and $M$ a $S$-module, a {\em commutative extension of 
$S$ by $M$} is an exact sequence:
$$
0\rightarrow M \stackrel{i}{\rightarrow} B \stackrel{p}{\rightarrow} S \rightarrow 0 \ ,
$$
$B$ being a commutative ring, $p$ a homomorphism of rings, and for all $b\in B$ and $m \in M$:
$bi(m)=i(p(b)m)$. Via $i$, $M$ is identified with an ideal in $B$, whose square is zero.
In fact, if $B$ is a ring and $I$ is a square zero ideal of $B$, then we can write the above
extension as
$$
0\rightarrow I \stackrel{i}{\rightarrow} B \stackrel{p}{\rightarrow} B/I \rightarrow 0 \ ,
$$
$I$ being an $R/I$-module. Abusively, we call also the ring $B$ {\em a commutative
extension of $S$ by $M$}.

Fossum proved in \cite{Fo} the folowing:
\begin{thrm}
A commutative extension:
$$
\exseq{E}{B}{S}
$$
of a Cohen-Macaulay ring $S$ by a canonical $S$-module $E$ is a 
Gorenstein ring.
\end{thrm}
\qed
\begin{re}
The multiplicity of the extension is double the multiplicity of the ring $S$,
since the multiplicity of the canonical module equals the multiplicity of the ring. 
This motivate us to call an extension like above {\bf a Fossum} doubling.
\end{re}

We can construct embedded doublings. 
If we take $B:=R/J$ and $S:=R/I$ as  quotients of a (say, regular) ring $R$, then an extension 
as above can be writen:
$$
0\rightarrow I/J \stackrel{i}{\rightarrow} R/J \stackrel{p}{\rightarrow} R/I \rightarrow 0 \ ,
$$
where $I^2 \subset J$.

The above explanations show that, given $S=R/I$ a Cohen-Macaulay ring, quotient of
the ring $R$, and $\pi :I/I^2\rightarrow E_S$ a surjection to a canonical $S$-module,
a commutative extension of $S$ by $E_S$ can be produced taking the kernel of $\pi$:
$$
0\rightarrow J/I^2 \rightarrow I/I^2 \stackrel{\pi}{\rightarrow} E_S \rightarrow 0 \ .
$$
The doubling of Fossum can be applied to schemes, {\em abstractly}, to obtain the so-called
{\bf ribbons} or {\em embedded}, to obtained the so-called Ferrand's doubling.
It is beyond the aim of this paper to speak about ribbons. To make the reader curious, we quote
from \cite{BE}: "The theory of ribbons is in some respects remarcably close to that of smooth
curves, but ribbons are much easier to construct and work with."
\subsection{ The geometric doubling (Ferrand)}

Consider now $X$ a locally Cohen-Macaulay subscheme of a regular one $Z$.
Then there is a dualizing $\O_X$-module $\omega _X$ which gives a good duality
on the category of coherent sheaves on $X$. Locally, $\omega _X$ corresponds to the canonical
module of the respective local rings of $X$. The method of Fossum can be then globalized
as follows:
 
\begin{constr}\label{Ferrand}
Let $X$ be a Cohen-Macaulay subscheme of the regular scheme $Z$. Let $I$ be the 
ideal (sheaf ideal, of course) defining $X$ in $Z$, 
$L$  a line bundle (locally free sheaf of rank 1) on $X$ and $p:I/I^2 \rightarrow 
\omega _X \otimes L$ a surjection. Then the kernel of $p$ is the of the form 
$J/I^2$, where $J$ is the ideal of a Gorenstein scheme
structure $Y$ in which $X$ is a closed subscheme. As sets, $Y=X$.
\end{constr}
This is the {\em Ferrand's doubling method}, \cite{Fe}.
We can write the exact sequence defining $J$ also:
$$
\exseq{L\otimes \omega _X}{\O_Y}{\O_X } \ .
$$
This exact sequence is the result of {\em patching together} exact sequences of the
Fossum's type.

In fact, Ferrand considers the case:

\begin{thrm}
In the above construction $Z=\P^3$, $L$ induced, i.e. $L=\O_X(n):=\O_{P^3}(n)|_X$.
If $p:I\rightarrow \omega _X(n)$ is a surjection such that $H^1(\P^3,I(-n))\rightarrow 
H^1(X,\omega _X)$ vanishes, then the curve $Y$ defined by $J=\hbox{ker}(p)$ is the scheme of
zeros of a rank $2$ vector bundle, $\omega _Y \cong \O_Y(-n)$ and $Ann(I/J)=I$. 
\end{thrm}

It is interesting to quote Ferrand (cf. \cite{Fe}): "L'id\'ee d'imposer la condition 
$Ann(I/J)=I$ provient de la th\'eorie de la liaison de Peskine et Szpiro, \cite{PS}". 
(For further connections with the linkage theory see the next sections.)

He combines the above construction with the following known facts:

1) {\em a locally Cohen-Macaulay curve $Y$ in $\P^3$ is the scheme of zeros of a rank $2$
vector bundle $E$ iff $\omega _Y \cong \O _Y(m)$ for some $m\in \Z$. $Y$ is globally complete
intersection iff $E$ splits} (cf. \cite{Se2}.

2) {\em If $E$ is a vector bundle of rank at least $2$ on a curve $X$, generated by its global 
sections, then there is a section $\O_X \rightarrow E$ which vanishes nowhere}.

to obtain

\begin{thrm}
If $X$ is a locally complete intersection curve in $\P^3$, defined by the ideal $I$,
then there exists a curve $Y$, containing $X$, which is the scheme of zeros of a section of
a rank $2$ vector bundle $E$ on $\P^3$. The ideal of $Y$ satisfies $Ann(I/J)=I$.
\end{thrm}

If we apply the above result for the affine space $Z=\A^3$, (cf. \cite{Sz}), as any vector 
bundle on $\A^3$ is trivial (cf. \cite{Su}, \cite{Q})
one answers affirmatively a question of Kronecker and Severi:

\begin{thrm}(Ferrand-Szpiro)(cf. \cite{Sz})
Any locally complete intersection curve in $\A^3$ is set-theoretically (i.e.
changing conveniently the scheme structure) a complete intersection.
\end{thrm}

\begin{re} 
The multiplicity of the local ring of $Y$ in a point is 
double of that of the same point considered on $X$. In particular, when $X$ is smooth the 
scheme $Y$ has multiplicity $2$ in any point.\end{re}

\section{Higher multiplicities}

The first description of the locally complete intersection scheme structures on a line
in $\P^3$ with multiplicity $3$ was given by Hartshorne, \cite{H1}. 

The first systematic study of multiple structures after the
papers of Fossum, Ferrand was done by B\u anic\u a \c si Forster in \cite{BF1}.
Unfortunately this paper had a small circulation. Here the classification up to multiplicity
$4$ is done, for scheme structures $Y$ on a smooth curve $X$ as suport, embedded in a smooth 
threefold $Z$.
The ideea was to cut (locally) the multiple structure with a transverse smooth germ of a 
surface. One obtains a multiple point in $\C\{x,y\}$. These are classified (cf. \cite{Br}. 
Deforming such a multiple point, one gets {\em a germ of a multiple curve}. The local results 
were then patched together to give complete constructions of the global multiple structures.
The main results of \cite{BF1} were included in \cite{BF2}, where a new ideea was introduced.
Namely, here the Cohen-Macaulay stratification, which was present also in \cite{BF1},
was constructed directly: if we cut the {\bf thick} scheme structure $Y$ on $X$ with the 
succesive infinitesimal neighbourhoods $X^{(i)}$ of $X$ (defined by the powers $I^i$ of the 
ideal $I$  of $X$ in $Z$) and {\bf throw away the embedded points}, (in dimension $1$, 
Cohen-Macaulay property is equivalent with the lack of embedded points), one gets a canonical 
filtration of $Y$ with Cohen-Macaulay curves. The possible filtrations are classified and one 
constructs conversely, step by step, from the reduced structure $X$  the thicker ones.

Independently (of \cite{BF2}), in \cite{M2} another Cohen-Macaulay
stratification was proposed, inspired from the linkage theory of Peskine and Szpiro \cite{PS}.
 Namely, if $J$ is the ideal of the
thick scheme $Y$ with $X$ as support, one considers the schemes defined by $J:I^i$.
Observe that this is not a genuine linkage, because in general $J\not\subset I^i$.
This was shown to work for curves and even more, for any dimension of $X$ and any dimension 
of $Z$, in the case of locally complete intersection multiple structure on smooth support,
up to multiplicity $4$.
The locally complete intersection (lci, for short) scheme structures {\bf in any dimension and 
any codimension, any characteristic of the base field} were classified, up to multiplicity $4$.

The two filtrations, from \cite{BF2} and \cite{M2} differ in general. For both of them
the succesive quotients of the defining ideals are, in good circumstances (e.g. {\em
on smooth support})
 {\bf vector bundles} on $X$. 
So, in principle if one takes succesively surjections $p_1:I\rightarrow E_1$, 
$p_2:\hbox{ker}(p_1)\rightarrow 
E_2$,\ldots , $p_m:\hbox{ker}(p_{m-1})\rightarrow E_m$, where $E_j$ are arbitrary
vector bundles, one gets all 
the thick structures in the asserted range, and a Cohen-Macauly filtration.
The problem is that such a procedure is by no means
canonical and the counting of parameters is difficult. On the contrary, the filtrations 
considered in \cite{BF2}
and \cite{M2} are canonical.
The one coming from the linkage theory has the advantage that provides directly some exact 
sequences. Also it can be applied to nonirreducible schemes. On the other 
hand the B\u anic\u a-Forster filtration has the advantage of a multiplicative structure. 
Namely, if $J_i:=$ the ideal of the the scheme obtained throwing away the embedded points of 
$Y\cap X^{(i)}$, then the graded object (with finitely many graded components) 
$\bigoplus J_i/J_{i+1}$ is a graded $\O_X$-algebra,
each piece of which is a vector bundle. 

The next step was done in \cite{M5}, where one shows that the filtration from \cite{M2}
works in fact up to multiplicity $6$ ({\bf any dimension, any codimension}) for lci
structures. In \cite{M5} a new filtration was introduced, namely $J:(J:I^i)$.
The associated graded object is a graded $\O_X$-algebra, $\Ac(Y)$. The graded object $\Mc(Y)$
associated to the filtration $J:I^i$ has a natural structure of graded $\Ac(Y)$-module. 
We denote by $\Bc(Y)$ the graded algebra coming from the B\u anic\u a-Forster filtration.
In general $\Bc(Y)$, $\Ac(Y)$, $\Mc(Y)$ are different. The multiplicative structure is always
nontrivial (of course, if the graded components involved are not zero). 
There are canonical homomorphisms $\Bc(Y)\rightarrow \Ac(Y)\rightarrow \Mc(Y)$.
In the rest of this section we discuss only the simplest constructions from \cite{M5}.
When one of the above graded objects has line bundles as components, then all three
filtrations are equal. Such a structure was called {\em quasiprimitive} in \cite{BF1}, 
\cite{BF2}. One shows in this case that, in general: $\Bc(Y)\cong\Ac(Y)\cong \Mc(Y)\cong
\oplus _{i=1}^{i=\mu -1}L^i(D_i)$, where $D_i$ are effective divisors on $X$ and $L$ is a
line bundle on $X$.
B\u anic\u a-Forster called a multiple structure with all $D_i=0$  {\em a primitive one}.
The fact that the multiplicative structure is not trivial implies inequalities 
fulfilled by $D_i$. The requirement that the structure is locally complete intersection 
restricts even more the shape of $\Ac(Y)$. The primitive structures are always lci 
(locally complete intersection). For the multiplicity $\mu =3$ no other quasiprimitive 
structure is lci. For $\mu =4$ there is only one possible shape of $\Ac(Y)$, besides the
primitives one: $\Ac(Y)=\O_X\oplus L \oplus L^2(D)\oplus L^3(D)$. For $\mu = 5$ there is
also only one possibility: $\Ac(Y)=\O_X\oplus L \oplus L^2(D)\oplus L^3(2D)\oplus L^4(2D)$.
For multiplicity $6$ the situation is more complicated: there are three models:
$$\Ac(Y)=\O_X\oplus L \oplus L^2\oplus L^3(D)\oplus L^4(D)\oplus L^5(D) \ ,$$
$$\Ac(Y)=\O_X\oplus L \oplus L^2(E)\oplus L^3(E)\oplus L^4(2E)\oplus L^5(2E) \ ,$$
$$\Ac(Y)=\O_X\oplus L \oplus L^2(F_1)\oplus L^3(F_1+F_2)\oplus L^4(2F_1+F_2)\oplus L^5(2F_1+F_2)
 \ .$$
which can be {\em combined}, in the following sense:
Consider $D$, $E$, $F_1$, $F_2$ effective divisors on $X$, such that $D$, $E$ and $F_1+F_2$
are pairwise disjoint. Then outside the support of $E$ and $F_1+F_2$ the shape is of the 
{\em type $D$}, outside $D$ and $F_1+F_2$ is of {\em type $E$}, etc. An example of such a
combination is:
$$\Ac(Y)=\O_X\oplus L \oplus L^2(E)\oplus L^3(D+E)\oplus L^4(D+2E)\oplus L^5(D+2E) \ $$
It is beyond the aim of this survey to present in detail the results of \cite{BF1}, \cite{BF2}
or \cite{M2}, \cite{M5}.

\section{Old and New Examples}

\subsection{Smooth Hilbert scheme in a point corresponding to double structures}
In \cite{BM} we proved the smoothness of the moduli space of 
stable rank $2$ vector bundles on $\P ^3$ with Chern classes $c_1=-1$, $c_2=4$,
in the points corrsponding to those vector bundles which correspond to double structures on
conics. A step in this proof was the fact that the Hilbert scheme $\H_{4t+6}$ of closed
subschemes of $\P^3$ with Hilbert polynomial $4t+6$ is smooth in the points corresponding
to double conics. {\bf The same proof}, which we repeat here {\em mutatis mutandi}, 
in the longer variant from \cite{BM'}, works for the following general case:
\begin{thrm} The Hilbert scheme $\H_{4t+r+2}$, of closed subschemes in $\P^3$ with Hilbert 
polynomial $4t+r+2$, $r> 0$, is smooth in the points corresponding to double conics and 
the set
$\H_d$ corresponding to double conics is open in $\H_{4t+r+2}$. The closure of $\H_d$ is an
irreudcible component of the Hilbert scheme.
\end{thrm}
\Proof  Let $C$ be a smooth conic in $\P^3$. Let $Y$ be a doubling given by an extension:
$$
\exseq{L}{ \O_Y} {\O_C}\ ,
$$
or, equivalently, by an exact sequence:
$$
\exseq{I_Y/I_C^2}{I_C/I_C^2}{L}\ ,
$$
where $L$ is a line bundle on $C$, which restricted to $\P^1$ via 
$\P^1\stackrel{i}{\cong} C\mono \P^3$ is $\O_{\P^1}(r)$.
As $i^*\O_C(\ell)=\O_{\P^1}(2\ell)$, the Hilbert polynomial of $Y$ is $\chi (\O_Y(t))=
\chi (\O_C(t))+\chi (L)=(2t+1)+(2t+r+1)=4t+r+2$. As $H^1(C,L)=H^1(\P^1,\O(r))=0$, one has
$\Pic (Y)\cong \Pic (C)$. Note that $\omega _Y|_C \cong \omega _C\otimes L^{-1}$ and then
$\omega _Y|_{\P^1}\cong \O _{\P^1}(-r-2)$.

To double  the conic $C$ "with $L$" is equivalent to giving 
$a\in H^0(L(1))\cong H^0(\O_{\P^1}(r+2)$ and $b\in H^0(L(2))\cong H^0(\O_{\P^1}(r+4)$
without common zeros on $C$, the pair $(a,b)$ being unique for each $Y$ with support $C$, up
to a factor in $\C^*$. (When $r$ is even, $r=2R$, is very easy to write the equations of $Y$.
The ideal of $Y$ is of the shape $(Ah+B\tau ,h^2,h\tau ,\tau ^2)$, where
$A\in H^0(L(1))=H^0(\O_C(R+1)$, $B\in H^0(L(2))=H^0(\O_C(R+2)$ have no common zeros along
$C$ and where $(h,\tau )$ are the equations which give $C$ in $\P^3$.) Consider the family
$\H^*=\{H,C,(a,b)\}/ \C^*$, where $H$ is a plane in $\P^3$, $C$ a conic in it and $(a,b)$ are
the above elements, identified with homogeneous polynomial in two variables, of degrees
$r+4$, respectively $r+3$, without common factors, modulo a factor in $\C^*$.
 $\H^*$ has a natural structure of an algebraic variety, which is connected, smooth,
quasiprojective, rational, of dimension $2r+15$ ($8$ parameters for the conic and $2r+7$ for the
doubling). The map $\H^* \rightarrow \H_{4t+r+2}$ 
which associate to $(H,C,(a,b))$ the doubling of $C$ with data $(a,b)$ is clearly algebraic.
It is easy to see that it is also bijective onto the locus $\H_d \subset \H_{4t+r+2}$ of double 
conics (cf. \cite{BM'}, {\em Lemma 5}). To show that $\H^*\cong \H_d$
we have to show that the tangent space at any point $y$ corresponding to double structure $Y$
on a conic $C$ has dimension $2r+15$. But, if $\H$ denotes our Hilbert scheme, $\dim T_{\H,y}=
h^0(N_Y|_{\P^3})=h^0((I_Y/I_Y^2)^{\vee})=h^1((I_Y/I_Y^2)\otimes \omega _Y)$
We shall use the exact sequences:
\begin{equation}
\exseq{(I_CI_Y/I_Y^2)\otimes \omega _Y}{(I_Y/I_Y^2)\otimes \omega _Y}
{(I_Y/I_CI_Y)\otimes \omega _Y}\label{1}
\end{equation}
\begin{equation}
\exseq{((I_Y^2+I_C^3)/I_Y^2)\otimes \omega _Y}{(I_CI_Y/I_Y^2)\otimes \omega _Y}
{(I_CI_Y/(I_Y^2+I_C^3))\otimes \omega _Y}\label{2}
\end{equation}
\begin{equation}
\exseq{(I_C^2/I_CI_Y)\otimes \omega _Y}{(I_Y/I_CI_Y)\otimes \omega _Y}
{(I_Y/I_C^2)\otimes \omega _Y}\label{3}
\end{equation}
From the exact sequence:
$$
\exseq{I_Y/I_C^2}{I_C/I_C^2\ (=\O_C(-1)\oplus \O_C(-2))}{I_C/I_Y\ (=L)}
$$
one deduces 
$$
I_Y/I_C^2\cong L^{-1} (-3)\ .
$$
This together with $I_C^2/I_CI_Y\cong (I_C/I_Y)^{\otimes 2}=L^2$ and $\omega _Y|_C\cong
L^{-1}(-1)$ allow us to rewrite (\ref{3}):
\begin{equation}
\exseq{L^1(-1)}{(I_Y/I_CI_Y)\otimes \omega _Y}{L^{-2}(-4)}\label{3a}
\end{equation}
In (\ref{2}) we need to compute the middle term. It is easy to see that the 
left term is given
by
$I_Y^2+I_C^3/I_Y^2\cong I_C^3/I_C^2I_Y\cong (I_C/I_Y)^{\otimes 3}=L^3$, and so:
$$
I_Y^2+I_C^3/I_Y^2\otimes \omega _Y \cong L^3\otimes \omega _Y \cong L^3\otimes \omega _Y|_C
\cong L^2(-1)\ .
$$
 For the right term we use the exact sequence:
$$
\exseq{((I_Y^2+I_C^3)/I_C^3)}{(I_CI_Y/I_C^3)}
{(I_CI_Y/(I_Y^2+I_C^3))}\ ,
$$
where the first term is:
$(I_Y^2+I_C^3)/I_C^3\cong I_Y^2/I_Y^2\cap I_C^3\cong I_Y^2/I_YI_C^2\cong 
(I_Y/I_C^2)^{\otimes 2}
\cong L^{-2}(-6)$,
and were we need also $\det (I_CI_Y/I_C^3)$. The exact sequence:
$$
\exseq{I_CI_Y/I_C^3}{I_C^2/I_C^3(\cong \O_C(-2)\oplus\O_C(-3)\oplus\O_C(-4))}
{I_C^2/I_CI_Y
(\cong L^2)}
$$
implies $\det (I_CI_Y/I_C^3)\cong L^{-2}(-9)$, so that $I_CI_Y/(I_Y^2+I_C^3)
\cong \O_C(-3)$, and then
$$
(I_CI_Y/(I_Y^2+I_C^3))\otimes \omega _Y \cong \O_C(-3)\otimes \omega _Y\cong L^{-1}(-4)
$$

The exact sequence (\ref{2}) becomes:
\begin{equation}
\exseq{L^2(-1)}{(I_CI_Y/I_Y^2) \otimes \omega _Y}
{ L^{-1}(-4)}\label{2a}
\end{equation}

Applying the long exact sequences of cohomology to (\ref{3a}), noticing:
\begin{equation}
H^0(L^{-2}(-4))\cong H^0(\P^1,\O(-2r-8))=0 \ ,
\end{equation}
\begin{equation}
H^1(L^{-2}(-4))\cong H^1(\P^1,\O(-2r-8))\cong \C^{2r+7}\ ,
\end{equation}
\begin{equation}
H^0(L(-1))\cong H^0(\P^1,\O(r-2))=\C^{r-1} \ ,
\end{equation}
\begin{equation}
H^1(L(-1))\cong H^1(\P^1,\O(r-2))\cong 0\ ,
\end{equation}
one gets:
\begin{equation}
H^0((I_Y/I_CI_Y)\otimes \omega _Y)\cong \C^{r-1}\ ,
\end{equation}
\begin{equation}
H^1((I_Y/I_CI_Y)\otimes \omega _Y)\cong \C^{2r+7}\ .
\end{equation}
Now apply the long exact sequence of cohomology to (\ref{2a}), noticing:
\begin{equation}
H^0(L^{2}(-1))\cong H^0(\P^1,\O(2r-2))=\C^{2r-1} \ ,
\end{equation}
\begin{equation}
H^1(L^{2}(-1))\cong H^1(\P^1,\O(2r-2))=0 \ ,
\end{equation}
\begin{equation}
H^0(L^{-1}(-4))\cong H^0(\P^1,\O(-r-8))=0 \ ,
\end{equation}
\begin{equation}
H^1(L^{-1}(-4))\cong H^1(\P^1,\O(-r-8))\cong \C^{r+7}\ ,
\end{equation}
and get:
\begin{equation}
H^0((I_CI_Y/I_Y^2)\otimes \omega _Y)\cong \C^{2r-1}
\end{equation}
\begin{equation}
H^1((I_CI_Y/I_Y^2)\otimes \omega _Y)\cong \C^{r+7}
\end{equation}

The long exact sequence of cohomology applied to (\ref{1}) gives:
\begin{eqnarray}
0\rightarrow H^0((I_CI_Y/I_Y^2) \otimes \omega _Y)\rightarrow H^0((I_Y/I_Y^2)
\otimes \omega _Y)
\rightarrow H^0((I_Y/I_CI_Y)\otimes \omega _Y)\rightarrow \nonumber \\
 H^1((I_CI_Y/I_Y^2) \otimes \omega _Y)\rightarrow H^1((I_Y/I_Y^2)\otimes \omega _Y)
\rightarrow H^1((I_Y/I_CI_Y)\otimes \omega _Y)\rightarrow 0 \ ,\nonumber
\end{eqnarray}
and the above computations imply that $h^{1}((I_Y/I_Y^2)\otimes \omega _Y) =2r+15$ is
equivalent to the injectivity of the map
$$
H^0((I_Y/I_CI_Y)\otimes \omega _Y)\rightarrow H^1((I_CI_Y/I_Y^2)\otimes \omega _Y)\ .
$$

We have the following commutative diagram with exact rows:
{\scriptsize
$$
\begin{CD}
0 @>>> (I_CI_Y/I_Y^2)\otimes \omega  @>>> (I_Y/I_Y^2)\otimes \omega  @>>> (I_Y/I_CI_Y)\otimes 
\omega  @>>> 0 \\
@.  @VVV @VVV @| @. \\
0 @>>> (I_CI_Y/I_Y^2+I_C^3)\otimes \omega  @>>> (I_Y/I_Y^2+I_C^3)\otimes \omega  @>>> 
(I_Y/I_CI_Y)\otimes \omega  @>>> 0 \\
@. @| @AAA @AAA @. \\
0 @>>> (I_CI_Y/I_Y^2+I_C^3)\otimes \omega  @>>> (I_C^2/I_Y^2+I_C^3)\otimes \omega  @>>> 
(I_C^2/I_CI_Y)\otimes \omega  @>>> 0
\end{CD}
$$
}
which gives the commutative diagram of cohomology:
$$
\begin{CD}
H^0((I_Y/I_CI_Y)\otimes \omega @>>> H^1((I_CI_Y/I_Y^2)\otimes \omega)\\
@| @VV\wr V \\
H^0((I_Y/I_CI_Y)\otimes \omega) @>>> H^1((I_CI_Y/I_Y^2+I_C^3)\otimes \omega)\\
@AA\wr A @|  \\
H^0((I_C^2/I_CI_Y)\otimes \omega @>>> H^1((I_CI_Y/I_Y^2+I_C^3)\otimes \omega)
\end{CD}
$$
such that we have to show that the lower row is injective, or equivalently $H^0((I_C^2/I_Y^2+I_C^3)
\otimes \omega _Y)=0$.

The locally free $\O_C$-module $I_C^2/I_Y^2+I_C^3$ can be computed from the exact sequence of 
$\O_C$-modules:
{\scriptsize
$$
\exseq{(L^{-2}(-6)\cong )\ (I_Y^2+I_C^3)/I_C^3}{I_C^2/I_C^3\ (\cong \O_C(-4)\oplus \O_C(-3)
\oplus \O_C(-2))}{I_C^2/(I_Y^2+I_C^3)}\ ,
$$
}
which pulled-back on $\P^1$ is isomorphic to:
$$
0\rightarrow \O_{\P^1}(-2r-12)\stackrel{\alpha}{\rightarrow}\O_{\P^1}(-8)\oplus \O_{\P^1}(-6)\oplus
\O_{\P^1}(-4)\stackrel{\beta}{\rightarrow}\O_{\P^1}(r-4)\oplus \O_{\P^1}(r-2)
\rightarrow 0\ ,
$$
where $\alpha =\left( \begin{array}{c} 
                a^2 \\
                2ab \\
		b^2 
		\end{array}\right)$
		and $\beta =\left( \begin{array}{ccc}
		             2b & -a & 0 \\
			     0  & -b & 2a
			     \end{array} \right) $.
Then $i^*(I_C^2/(I_Y^2+I_C^3)\otimes \omega _Y) \cong \O_{\P^1}(-6)\oplus \O_{\P^1}(-4)$, 
and so $H^0((I_C^2/I_Y^2+I_C^3)\otimes \omega _Y)=H^0(\P^1,\O_{\P^1}(-6)\oplus \O_{\P^1}(-4))=0$.     
\qed			     
\begin{re}
(i) The fact that $\H_d$ is an irreducible component of the Hilbert scheme was reobtained
in \cite{NS}, where, more than this, 
{\bf the whole decomposition in irreducible components of the Hilbert scheme}
 is described and, {\em more difficult}, {\bf the connectedness of the Hilbert scheme
 of curves of degree $4$} is proved. The dimension of $\H_d$ follows also from \cite{HS}.

(ii) It is not true in general that double structures on plane curves in $\P^3$ give smooth points
of the respective Hilbert scheme. (cf. \cite{BM'}).

(iii) The case of ropes on lines is analyzed in \cite{NNS} and a smoothness result is obtained
(cf. loc. cit. for the exact statement). Although in \cite{NNS} the method is not the computing of
the global sections of the normal sheaf, some parallel to our old proof can be observed
(e.g the use of matrices similar to our $\alpha $, $\beta $).

\end{re}		
\subsection{A Double or a Triple Plane in $\P ^5$ is Never
the Scheme of Zeros of a Section of an Indecomposable Vector Bundle on $\P ^5$}

We mentioned that any smooth curve in $\P^3$ is the support of a locally complete 
structure of double degree which is the scheme of zeroes of a section of a convenient 
rank $2$ indecomposable vector bundle on $\P^3$. For surfaces in $\P^5$ this is no 
more true. We show:
\begin{thrm}
Let $\P^2 \cong X\subset \P^5$ a plane in $\P^5$. There is no indecomposable rank three
vector bundle $F$ on $\P^5$ with a section having as vanishing scheme a double structure 
on $X$.
\end{thrm}
\Proof Suppose such a doubling $Y$ exists. After \cite{M2} a double structure on
a plane is given by a (Fossum-)Ferrand construction. 
Then we should have an exact sequence:
$$
\exseq{\O_X(r)}{\O_Y}{\O_X}\ ,
$$
where $r\ge 0$ (for $r=-1$ the doubling is globally complete intersection, so that $F$
exists and splits).
If $Y$ is the scheme of zeros of a section $s\in H^0(F)$ in a rank $3$ vector bundle $F$,
and $E:=F^\vee$, then one has an exact sequence:
$$
0\rightarrow \Lambda ^3E \rightarrow \Lambda ^2E \rightarrow E \rightarrow \O_{\P^5}
\rightarrow \O_Y \rightarrow 0 \ .
$$
The doubling exists for each $r$ in our range because there are surjections
$$
I_X/I_X^2 \cong 3\O_X(-1)\stackrel{p}{\rightarrow} \O_X(r)\ ,
$$
such a surjection being determined by three forms of degree $r+1$ in the homogeneous
coordinates on $X=\P^2$, without common zeros. (We used tacitely the notation of the type
$3A$ for $A\oplus A\oplus A$.)

The plan of the proof is the following: we compute the Hilbert polynomial of $Y$
from the two above exact sequences, as a function in $r$ and as a function in
the Chern classes of $E$ (one uses Riemann-Roch-Hirzebruch Theorem, cf. \cite{Hi}). 
This will determine the Chern classes of $E$ with respect
to $r$ and then we see that $E$ does not satisfy the Schwarzenberger conditions 
(cf. {\em Appendix One} by R.~L.~E.~Schwarzenberger to \cite{Hi} ---
they come from the requirement that the Hilbert-Euler characteristic of $E$ takes
values in $\Z$).

Denote by $c_1$, $c_2$, $c_3$ the Chern classes of $E$. 
Then $\Lambda ^3E\cong \O_{\P^5}(c_1)$, $\Lambda ^2E$ has Chern classes :
$c_1'=3c_1$, $c_2'=c_1^2+c_2$, $c_3'=c_1c_2-c_3$. 
After a long computation, which we don't reproduce here, one gets the Hilbert
polynomial of $Y$:
$$
\chi _Y(t)=-\frac{c_3}{2}t^2-\frac{(c_1+6)c_3}{2}t+
\frac{(c_2-2c_1^2-18c_1-51)c_3}{2}
$$
On the other hand the defining extension of $\O_Y$ gives:
$$
\chi _Y(t)=\chi _{\P^2}(t)+\chi _{\P^2}(t+r)={t+2\choose 2}+{t+r+2 \choose 2}=
t^2+(r+3)t+\frac{r^2+3r+4}{2}
$$
Comparing the two formulas one gets:
\begin{eqnarray}
c_1 & = & r-3 \ ,\nonumber \\
c_2 & = & \frac{3r^2+9r+26}{2} \ ,\nonumber \\
c_3 & = & -2 \ .\nonumber
\end{eqnarray}
We show now that there is no value of $r$ for which the above Chern classes
can be attained. Namely, we write the Hilbert polynomial of $E$ in the basis
${t+i \choose i},\ i=1,\ldots ,5$:
\begin{equation}\begin{split}
\chi _E(t) & =3{t+5\choose 5}-(r-3){t+4 \choose 4}+(r^2+7r+10){t+3\choose 3}+\\
&  \frac{7r^3+30r^2+29r-54}{12}{t+2\choose 2}+\frac{r^4-24r^3-197r^2-560r-548}
{48}{t+1\choose 1}\\ & -\frac{19r^5+235r^4+1305r^3+3765r^2+5616r+3140}{480}
\end{split}\end{equation}
The polynomial $\chi _E(t)$ takes integer values for $t\in \Z$ iff the coefficients
are in $\Z$. In particular we should have:
\begin{eqnarray}
7r^3+30r^2+29r-54 \equiv 0 \ \ ({\rm mod}\ 3)\nonumber \\
r^4-24r^3-197r^2-560r-548 \equiv 0 \ \ ({\rm mod}\ 3)\nonumber
\end{eqnarray}
what is impossible.
\begin{re}
The above principially easy observation should be very well known, as the 
construction of vector bundles on the projective space is a much desired objective.
We don't know any reference for it.
\end{re}
\qed
\begin{thrm}
Let $\P^2 \cong X\subset \P^5$ a plane in $\P^5$. There is no indecomposable rank three
vector bundle $F$ on $\P^5$ with a section having as vanishing scheme a triple structure 
on $X$.
\end{thrm}
\Proof According to \cite{M2} a triple structure, which is locally
a complete intersection, supported by a smooth surface, is similar to a 
{\em triple primitive structure} on a curve (in the B\u anic\u a-Forster terminology , 
cf. \cite{BF1}, \cite{BF2}). That means that a triple structure $Z$ on $X$
contains a double one $Y$ and one has the defining exact sequences:
\begin{equation}
\exseq{I_Y/I_X^2}{I_X/I_X^2}{L\ (=\O_X(r))}\ \label{double},
\end{equation}
$$
0\rightarrow I_Z/I_XI_Y\rightarrow I_Y/I_XI_Y \stackrel{p}{\rightarrow} \O_X(2r)\rightarrow 0\ ,
$$ 
where $p$ is a retract of the natural inclusion $L^2=I_X^2/I_XI_Y\hookrightarrow I_Y/I_XI_Y$.
\begin{cl}
The above retract does exist.
\end{cl}
{\em Proof } of the claim. Denote by $G$ the rank $2$ vector bundle $I_Y/I_X^2$ on $X$.
The canonical exact sequence:
$$
\exseq{I_X^2/I_XI_Y}{I_Y/I_XI_Y}{I_Y/I_X^2}
$$
determines an element $\xi \in {\rm Ext}^1(G,L^2)$. But ${\rm Ext}^1(G,L^2)\cong 
H^1(G^\vee (2r))=0$, as one sees from the exact sequence (\ref{double}), and so the 
above exact sequence splits.
\qed

We come back to the proof of our theorem. The proof will be anlogous to the proof of the
previous result.

The Hilbert polynomial of $Z$ is:
$$
\chi _Z(t)=\frac{3}{2}t^2+\frac{6r+9}{2}t+\frac{5r^2+9r+6}{2}
$$
On the other hand the existence of an exact sequence:
$$
0\rightarrow \Lambda ^3E \rightarrow \Lambda ^2E \rightarrow E \rightarrow \O_{\P^5}
\rightarrow \O_Z \rightarrow 0 
$$
imply again:
$$
\chi _Z(t)=-\frac{c_3}{2}t^2-\frac{(c_1+6)c_3}{2}t+
\frac{(c_2-2c_1^2-18c_1-51)c_3}{2} \ ,
$$
where $c_i$ are the Chern classes of $E$. Compairing the two results, one gets:
\begin{eqnarray}
c_1 & = & 2r-3 \nonumber \ ,\\
c_2 & = & \frac{19r^2+27r+39}{3}\nonumber \ ,\\
c_3 & = & -3 \ .\nonumber
\end{eqnarray}
The expression of $c_2$ shows that $r$ is a multiple of $3$. Denote $r=3R$.
Then:
\begin{eqnarray}
c_1 & = & 6R-3 \nonumber \ ,\\
c_2 & = & 57R^2+27R+13 \nonumber \ ,\\
c_3 & = & -3 \ .\nonumber
\end{eqnarray}
This shape of the Chern classes is impossible, because the Hilbert polynomial of 
$E$ computed by the Riemann-Roch-Hirzebruch Theorem does not take integral values.
Namely, its coefficient of ${t+1\choose 1}$ in the developing in the basis
${t+i\choose i},\  i=1,\ldots 5$, which is
$$
\frac{207R^4-1512R^3-1845R^2-828R-134}{12}
$$
is not in $\Z$.
\qed

\noindent Nicolae Manolache\\
Institute of Mathematics "Simion Stoilow"\\
of the Romanian Academy \\
P.O.Box 1-764
Bucharest, RO-014700

\noindent e-mail: nicolae.manolache@imar.ro
\end{document}